\documentclass[12pt]{article}

\usepackage{amssymb,amsthm,amsmath}
\usepackage{url}

\def\eps{\varepsilon}

\def\E{\mathbb E}
\def\Ex{\mathbb E}

\def\PP{\mathbb P}
\def\P{\mathbb P}
\def\Pr{\mathbb P}
\def\er{\mathbb R}

\newcommand{\R}{\mathbb{R}}

\newcommand{\dist}{{\rm dist\,}}
\newcommand{\supp}{{\rm supp\,}}

\newtheorem{lemma}{Lemma}[section]

\newtheorem{theorem}[lemma]{Theorem}

\newtheorem{cor}[lemma]{Corollary}

\def\r{\right}

\def\no{\|\cdot\|}
\def\la{\langle}
\def\ra{\rangle}

\def\Gam{\Gamma}
\def\gam{\gamma}
\def \sig {\sigma}

\def\Id{\mathrm{Id}}
\def\conv{\mathrm{conv}}
\def\diam{\mathrm{diam}}

\date{}

\title{Chevet type inequality and norms of submatrices}

\author{{Rados{\l}aw Adamczak${}^{1}$, Rafa{\l} Lata{\l}a${}^{1}$, 
Alexander E. Litvak${}^{2}$,}\\
{Alain Pajor${}^{3}$, Nicole  Tomczak-Jaegermann${}^{4}$}}

\newcommand\address{\noindent\leavevmode
\noindent
Rados{\l}aw  Adamczak, \\
Institute of Mathematics, \\
University of Warsaw, \\
Banacha 2, 02-097 Warszawa, Poland\\ 
\texttt{\small e-mail:  radamcz@mimuw.edu.pl}

\medskip
\noindent
Rafa{\l} Lata{\l}a, \\
Institute of Mathematics, \\
University of Warsaw, \\
Banacha 2, 02-097 Warszawa, Poland\\ 
and \\
Institute of Mathematics,\\
Polish Academy of Sciences,\\
\'{S}niadeckich 8, 00-956 Warszawa, Poland\\
\texttt{\small e-mail:  rlatala@mimuw.edu.pl}

\medskip
\noindent
Alexander E. Litvak, \\
Dept.~of Math.~and Stat.~Sciences,\\
University of Alberta, \\
Edmonton, Alberta, Canada, T6G 2G1.\\
\texttt{\small e-mail:  alexandr@math.ualberta.ca}

\medskip
\noindent
Alain  Pajor, \\
Universit\'{e} Paris-Est\\
\'{E}quipe d'Analyse et Math\'{e}matiques Appliqu\'ees, \\
5, boulevard Descartes,
Champs sur Marne,\\
77454 Marne-la-Vall\'{e}e,  Cedex 2, France\\
\texttt{\small e-mail: Alain.Pajor@univ-mlv.fr }

\medskip
\noindent
Nicole  Tomczak-Jaegermann, \\
Dept.~of Math.~and Stat.~Sciences,\\
University of Alberta, \\
Edmonton, Alberta, Canada, T6G 2G1.\\
\texttt{\small e-mail:    nicole.tomczak@ualberta.ca}

}

\begin{document}

\maketitle
\footnotetext[1]{Research partially supported by MNiSW Grant no. N N201 397437 
and the Foundation for Polish Science.}
\footnotetext[2]{Research partially supported by  the 
E.W.R. Steacie Memorial Fellowship.}
\footnotetext[3]{Research partially supported by the ANR project ANR-08-BLAN-0311-01.}
\footnotetext[4]{This author holds the Canada Research Chair in
  Geometric Analysis.}

\begin{abstract}  
We prove a Chevet type inequality which gives an upper bound for the 
norm of an isotropic log-concave unconditional random matrix in terms 
of  expectation of the supremum of ``symmetric exponential" processes
compared to the Gaussian ones in the Chevet inequality. This is used to give 
sharp upper estimate for a quantity $\Gam_{k,m}$ that controls uniformly the 
Euclidean operator norm of the sub-matrices with $k$ rows and $m$ columns of 
an isotropic log-concave unconditional random matrix. We apply these estimates 
to give a sharp bound for the Restricted Isometry Constant of a random matrix
with independent log-concave unconditional rows. We show also that our Chevet 
type inequality does not extend to general isotropic log-concave random matrices. 
\end{abstract}

\section{Introduction}
\label{intro}

Let $n$, $N$ be positive integers. 
Let $K\subset \mathbb{R}^N$ and  $L\subset \mathbb{R}^n$
be origin symmetric convex bodies,  $\| \cdot \|_K$ and $\| \cdot \|_L$
be the corresponding gauges  on $\mathbb{R}^N$ and $\mathbb{R}^n$,
that is the norms for which $K$ and $L$ are the unit balls.

To shorten the notation we write $\|\Gam: K\to L\|$ for the operator norm of 
a linear operator $\Gam: (\R^N, \|\cdot\|_K)\to (\R^n, \|\cdot\|_L)$. In particular, 
$\|\Gam : K\to B_2^N\|$ will denote the operator norm of $\Gam$ considered as a 
linear operator from $(\R^N, \|\cdot\|_K)$ to $\ell_2^N$, where $\ell_2^N$ 
is $\R^N$ equipped with the canonical Euclidean norm, whose unit ball is $B_2^N$; 
similarly for $\|\Gam : B_2^n\to L\|$. Note also that the dual normed space 
$(\R^N, \|\cdot\|_K)^*$ of $(\R^N, \|\cdot\|_K)$ may be identified 
(via the canonical inner product) with $(\R^N, \|\cdot\|_{K^{\circ}})$,  
where ${K^\circ}$ denotes the polar of $K$ 
(see the next section for all definitions). The canonical basis on
$\mathbb{R}^d$  is denoted by $\{e_i\}_{1\leq i\leq d}$.

Let $(g_i)_{1\leq i\leq \max{(n,N)}}$ be i.i.d. standard Gaussian 
random variables that is centered Gaussian variables with variance 1,
and $\Gam$ be a Gaussian matrix whose entries are i.i.d. standard Gaussian. 
Then one side of the Chevet inequality (\cite{Chev}, see also \cite{G}  
for sharper constants) states that
\begin{align}
\mathbb{E} \|\Gam : K \to L \|
 &\leq
  C \|\Id : K \to B_2^N \|
  \cdot \mathbb{E}\left\| \sum_{i=1}^n g_i e_i \right\| _{L} \nonumber
\\   &+
  C \|\Id : B_2^n \to L \|  \cdot
   \mathbb{E}\left\| \sum_{i=1}^N g_i e_i \right\| _{K^\circ} ,
\label{chevet1}
\end{align}
where $\Id$ stays for the formal identity operator  
and $C$ is an absolute constant. This 
inequality plays an important role in Probability in Banach Spaces and 
in Asymptotic Geometric Analysis (\cite{BG, Tom}).

We say that a random matrix $\Gam=(\gamma_{ij})$ is {\em isotropic} if 
all entries $(\gamma_{ij})$ are uncorrelated centered with variance one 
and it is {\em log-concave} if the joint distribution of the $\gamma_{ij}$'s 
has a density which is log-concave on its support, finally we say that the 
matrix $\Gam$ is {\em unconditional} if  for any 
choice of signs $(\varepsilon_{ij})$ the matrices $\Gam$ and
$(\varepsilon_{ij} \gamma_{ij})$ have the same distribution. 
There are similar definitions for random vectors.

 In Theorem~\ref{Chevet} we prove that an inequality similar to the 
Chevet inequality (\ref{chevet1}) holds for any isotropic 
log-concave unconditional random matrix $\Gam$
when substituting the Gaussian random variables $g_i$'s
by i.i.d. random variables with symmetric
exponential distribution with variance 1. Moreover, in 
Corollary~\ref{probest} we provide the corresponding probability 
estimates.

A result from \cite{La2} of the second named author
of this article states that if
$X=(X_1,\dots, X_d)$ is an isotropic log-concave unconditional random vector in
$\R^d$
and if $Y = (E_1, \ldots, E_d)$,
where $E_1, \ldots, E_d$ are i.i.d. symmetric exponential random variables,
then for any norm $\|\cdot\|$ on $\R^d$, one has

\begin{equation}\label{latala1}
     \E\|X\|  \le C\ \E\|Y\|,
\end{equation}
where $C$ is an absolute constant.

The proof of our Chevet type inequality consists of two steps.
First, using the comparison (\ref{latala1}),
we reduce the case of a general isotropic log-concave unconditional
random matrix $A$  to the case of an exponential  random
matrix, i.e. the matrix whose entries are i.i.d. standard symmetric
exponential random variables.
The second step uses Talagrand's result (\cite{TalCan}) on relations
between some random processes associated to the symmetric exponential distribution
and so-called $\gamma_p$ functionals.

We apply our inequality of Chevet type to obtain sharp uniform bounds
on norms of sub-matrices of isotropic log-concave unconditional
random matrices $\Gam$. More precisely,
for any subsets $J\subset \{1,\ldots,n\}$ and $I \subset \{1,\ldots,N\}$
denote the submatrix of $\Gam$ consisting of the rows indexed by elements
from $J$ and the columns indexed by elements from $I$ by $\Gam(J, I)$.
Given $k\le n$ and $m\le N$ define the parameter $\Gam_{k,m}$ by
$$
 \Gam_{k,m} =  \sup  \|\Gam (J, I)\ :\ \ell_2^m \to \ell_2^k\| ,
$$
where the supremum is taken over all subsets
$J\subset \{1,\ldots,n\}$ and $I \subset \{1,\ldots,N\}$  with
cardinalities  $|J| = k$, $|I|= m$. That is, $\Gam_{k,m}$ is
the maximal operator norm of a sub-matrix of $\Gam$ with $k$ rows and
$m$ columns.
We prove that
$$
   \Gam_{k,m} \leq C\left(  \sqrt{m}\log\left(\frac{3 N}{m}\right)
  + \sqrt{k}\log\left(\frac{3 n}{k}\right)\r),
$$
with high probability. This estimate is sharp up to absolute constants.

Furthermore, we provide applications of this result to
the Restricted Isometry Property (RIP) of a matrix with
independent isotropic log-concave unconditional random rows.
We give sharp estimate for the restricted isometry constant
of such matrices.

It is well known and follows from Talagrand's majorizing
measure theorem (see \cite{Tal}) that if
$X=(X_1,\dots, X_d)$ is a centered sub-gaussian random vector in
$\R^d$ with parameter $\alpha>0$, that is, all coordinates $X_i$ 
are centered and  for any $x\in\R^d$ of Euclidean norm 1, any $t>0$,
$\P(|\sum x_iX_i|\geq t)\le 2 \exp(-t^2/\alpha^2)$, then for any norm 
$\|\cdot\|$ on $\R^d$, one has
\begin{equation}\label{talagrand1}
     \E\|X\|  \le C\alpha\, \E\|Y\|,
\end{equation}
where $Y = (g_1, \ldots, g_d)$ 
and $C>0$ is an absolute constant.

It is interesting to view both inequalities
(\ref{latala1}) and (\ref{talagrand1}) in parallel. There are both based on
majorizing measure theorems of Talagrand; inequality (\ref{talagrand1}) states that
the expectation of the norm of a sub-gaussian vector is
up to a multiplicative constant, dominated by its Gaussian replica.
So Gaussian vectors are almost maximizers. To which class of random vectors does 
inequality (\ref{latala1}) correspond? 
Since in many geometric and probabilistic inequalities involving isotropic 
log-concave vectors, Gaussian and exponential vectors are the extreme cases, 
it was naturally conjectured that the expectation of the norm of isotropic 
log-concave vector is similarly dominated by the corresponding expectation of 
the norm of an exponential random vector. This conjecture would have many 
applications. For instance the estimate of $\Gam_{k,m}$ above would extend to general 
log-concave random matrices, which is open (see  \cite{ALLPT3}). 


We show that this is not the case. Namely, in Theorem~\ref{latexam} we prove that  
for any $d\geq 1$, there exists an isotropic log-concave random vector $X\in\R^d$ 
and a norm $\|\cdot\|$ on $\R^d$ such that
\begin{equation}\label{latala2}
     \E\|X\|  \ge c  \sqrt{\ln d}\, \E\|Y\|, 
\end{equation}
where $Y$ is of ``symmetric exponential"  type and $c$ is a positive universal 
constant. Similarly we show that our Chevet inequality does not extend to the 
setting of general log-concave random matrices (non unconditional). In fact 
it would be interesting to find the best dependence on the dimension in the 
reverse inequality to (\ref{latala2}). More precisely, to solve the following problem.

\medskip

\noindent
{\bf Problem. }{\it  Find tight (in terms of dimension $d$) estimates 
for the following quantity 
$$
   C(d) = \sup _{\no} \sup _X \frac{\E\|X\|}{\E\|Y\|},
$$
where $Y = (E_1, \ldots, E_d)$ and the supremum is taken over all norms 
$\no$ on $\R^d$ and all isotropic log-concave random vectors $X\in\R^d$.
}

\medskip

Theorem~\ref{latexam} and Remark~2 following it show  that 
$
   c\ \sqrt{\ln d} \leq C(d) \le C \sqrt{d} 
$
 for some absolute positive constants $c$ and $C$.

The results on norms of submatrices and applications were partially announced in \cite{ALLPT4}.
For the related estimates in the non-unconditional case,
see \cite{ALLPT3}.

The paper is organized as follows. In the next section we introduce notation and quote 
known results which will be used in the sequel. In Section~\ref{chevetineq} we prove 
the Chevet type inequality (and corresponding probability estimates) for unconditional 
log-concave matrices. In remarks we discuss its sharpness showing that in general one 
can't expect the lower bound of the same order and providing a relevant lower bound. 
In Section~\ref{RIP} we apply our Chevet type inequality to obtain  sharp uniform 
estimates for norms of submatrices. Then we apply the results to the RIP. 
Section~\ref{example} is devoted to  examples showing that one can't drop the condition 
of unconditionality in the comparison theorem of the second named author and in our 
Chevet type inequality. Finally, in Section~\ref{dirproof}, we present a direct approach 
to uniform estimates of norms of submatrices, which does not involve Chevet type 
inequalities and $\gamma_p$ functionals, but is based  only on tail estimates for suprema 
of linear combinations of independent exponential variables and on a chaining argument 
in spirit of \cite{ALPT}.

\medskip

\noindent
{\bf Acknowledgment:\ } The research on this project was partially
done when the authors participated in the Thematic Program on
Asymptotic Geometric Analysis at the Fields Institute in Toronto in
Fall 2010 and in the Discrete Analysis Programme at the Isaac Newton
Institute in Cambridge in Spring 2011. The authors wish to thank these
institutions  for their  hospitality and excellent working conditions.

\section{Notation and Preliminaries}
\label{notat}

 By $|\cdot|$ and $\la \cdot , \cdot \ra$ we denote the canonical
Euclidean norm and the canonical inner product on $\R ^d$.
The canonical basis of $\R ^d$ is denoted by $e_1, \ldots, e_d$.

As usual, $\| \cdot \| _p$, $1\leq p \leq \infty$, denotes the
$\ell _p$-norm, i.e. for every $x=(x_i)_{i=1}^d \in\R^d$ 
$$
 \|x\| _p = \left( \sum _{i= 1}^d |x_i|^p \r) ^{1/p} \,
 \mbox{ for } \ p < \infty \, \, \, \, \mbox{ and } \, \, \, \,
 \|x\| _{\infty } = \sup _{i\leq d} |x_i|
$$
and $\ell _p^d = (\R ^d, \|\cdot \|_p)$.  
The unit ball of $\ell _p^d$ is denoted  by $B_p^d$. 
For a non-empty set $T\subset \R^d$ we write
$\diam_p(T)$ to denote the diameter of $T$ with respect to 
the $\ell_p$-norm.

For an origin symmetric convex body $K\subset \R^d$, the Minkowski 
functional of $K$ is
$$
    {\| x\|}_{K}=\inf \{\lambda >0 \ |\ x\in\lambda K\},
$$
i.e. the norm, whose unit ball is $K$. The polar of $K$ is 
$$
 K^{\circ} = \{x \ | \ \la  x, y \ra \leq 1 \ \ \mbox{ for all } \ y\in K \}. 
$$
Note that $K^{\circ}$ is the unit ball of the space dual to $(\R^d, \|\cdot\|_K)$.

Given an $n\times N$ matrix $\Gam$ and origin symmetric convex bodies $K\subset \R^N$,
$L\subset \R^n$ we denote by
$$
  \|\Gam : K \to L \|
$$
the operator norm of $\Gam $ from $(\R^N,  {\| \cdot\|}_{K})$ to $(\R^n,  {\| \cdot\|}_{L})$.
We also denote 
$$
   R(K) = \|\Id : K\to B _2^N \|, \quad   R(L^{\circ}) = \|\Id : B_2^n \to L \|
  = \| \Id : L^{\circ} \to B_2^n\|, 
$$
where $\Id$ denotes the formal identity $\R^N \to \R^N$ or  $\R^n \to \R^n$.

Given a subset $K\subset \R^d$ the convex hull of $K$ is denoted by $\conv (K)$.  

A random vector $X =(X_1,\ldots,X_N)$ is called
unconditional if for every sequence of signs
$\varepsilon_1,\ldots,\varepsilon_N$, the law of
$(\varepsilon_1 X_1,\ldots,\varepsilon_N X_N)$ is the same as the
law of $X$.

A random vector $X$ in $\R^n$  is called isotropic if
$$
\E\langle X,y\rangle=0,\quad \E\,|\langle X, y \rangle|^{2}=\|y\|_2^{2}
\quad \mbox{\rm for all }
 y\in \R^{n},
$$
in other words, if $X$ is centered and its covariance matrix
$\E\, X\otimes X$ is the identity.

A random vector $X$  in $\R^n$ with full dimensional support is called 
log-concave if it has a log-concave density. Notice that all isotropic 
vectors have full dimensional support.


By $E_i$, $E_{ij}$ we denote independent symmetric exponential
random variables with variance $1$ (i.e. with the density 
$2^{-1/2} \exp(-\sqrt{2}\ |x|)$).  By $g_i$, $g_{ij}$ we denote 
standard independent ${\cal N}(0, 1)$ Gaussian random variables. 
The $n\times N$ random matrix with entries $g_{ij}$ will be called 
the Gaussian matrix, the $n\times N$ random matrix with 
entries $E_{ij}$ will be called the exponential random matrix.
Similarly, the vectors $G=(g_1, \ldots, g_d)$ and 
$Y=(E_1, \ldots, E_d)$ are called Gaussian and exponential 
random vectors.

In the sequel we often consider $n\times N$ matrices as 
vectors in $\R^d$ with $d=nN$ and the inner product defined by 
$$
  \la A, B\ra = \sum _{i, j} a_{i j} b_{i j} 
$$
for $A = (a_{i j})$, $B = (b_{i j})$. Clearly, 
the corresponding Euclidean structure is given by Hilbert-Schmidt 
norm of a matrix: 
$$
 |A| =  \|A\| _2 = \left(\sum _{i, j} |a_{i j}|^2 \r)^{1/2} . 
$$
In this notation we have $\|A\| _{\infty} = \max _{i, j} |a_{i j}| $. 
We say that such a matrix $A$ is isotropic/log-concave/unconditional if 
it is isotropic/log-concave/unconditional as a vector in $\R^d$, $d=nN$
(cf. the definition given in the introduction). 

Given $x\in \R^N$ and $y\in \R^n$, denote by $x\otimes y =y x^\top$ the matrix 
$\{y_i x_j\}_{ij}$, i.e. the matrix corresponding to the linear operator 
defined by
$$
  x\otimes y  \ (z) = \la z , x \ra  y.
$$ 
Then, for an $n\times N$ matrix $\Gam$ 
$ = (\gamma _{i j} )$,
$$
   \|\Gam : K \to L \| = \sup _{x\in K} \sup _{y\in L^{\circ}}
    \sum _{i, j} \gamma _{i j} x_j y_i = \sup _T
    \la \Gam, x\otimes y \ra , 
$$
where the latter supremum is taken over
$$
   T = K\otimes L^{\circ} = \{ x\otimes y \ \colon \ x\in K,\ y\in L^{\circ}  \}.
$$

We will use the letters $C, C_0, C_1, \ldots$, $c, c_0, c_1, \ldots$ 
to denote positive absolute constants whose values may differ at each 
occurrence. 
We also use the notation $F\approx G$ if there are two positive absolute 
constants $C$ and $c$ such that $c \, G \leq F \leq C \, G$. 

\medskip

Now we state some results which will be used in the sequel. 
We start with the following lemma, which provides asymptotically 
sharp bounds on the norm of the exponential matrix considered as 
an operator $\ell _1^N \to \ell _1^n$. We will use it in our 
examples on sharpness of some estimates.

\begin{lemma}\label{lonenorm}
Let 
$\Gamma =(E_{ij})_{i \leq n, j\leq N}$. Then
$$
  \E \ \| \Gam \ :\ \ell _1^N \to \ell _1^n \| 
  \approx n + \ln N.
$$
\end{lemma}


\noindent
{\bf Proof.} 
 First note 
\begin{equation}\label{llonenorm}
  \| \Gam \ :\ \ell _1^N \to \ell _1^n \|=\max_{i\leq N}\sum_{j=1}^n|E_{ij}|.   
\end{equation}
By the Chebyshev inequality for every $i\leq n$ we have  
$$
  \Pr\Big(\sum_{j=1}^n|E_{ij}|\geq t\Big) \leq \exp\Big(-\frac{t}{2}
  \Big)\ \Ex\exp \Big(\frac{1}{2}\sum_{j=1}^n|E_{ij}|\Big) \leq 
  C^n\exp\Big(-\frac{t}{2}\Big) 
$$
for some absolute constant $C>0$. 
Hence the union bound and integration by parts gives
$$
  \E  \ \| \Gam \ :\ \ell _1^N \to \ell _1^n \| 
  \leq  C\left(n + \ln N \r) .
$$

On the other hand,  by (\ref{llonenorm})
$$
  \Ex\| \Gam \ :\ \ell _1^N \to \ell _1^n \|\geq \Ex\sum_{j=1}^n|E_{1j}| = n/\sqrt{2}   
$$
and 
$$
  \Ex\| \Gam \ :\ \ell _1^N \to \ell _1^n \|\geq \Ex\max_{i\leq N} |E_{i1}| 
  \approx 1+ \ln N 
$$
(the last equivalence is well-known and follows from direct computations). 
This completes the proof. 
\qed

\medskip

The next theorem is a comparison theorem from \cite{La2}. 

\begin{theorem} \label{Latala}
Let $X$ be an isotropic log-concave unconditional random vector in
$\R^d$ and $Y = (E_1, \ldots, E_d)$ be an exponential random vector. 
Let $\|\cdot\|$ be a norm on $\R^d$. Then
$$
 \E\|X\|  \le C\ \E\|Y\| , 
$$
where $C$ is an absolute positive constant. 
Moreover, for every $t \ge 1$, 
$$
  \P(\|X\| \ge t) \le C\ \P(\|Y\|\ge t/C).
$$
\end{theorem}

\medskip

\noindent
{\bf Remark. }
The condition ``$X$ is unconditional" cannot be omitted  in 
Theorem~\ref{Latala}. We show  an example 
proving that in Section~\ref{example}. 

\medskip

We will also use two Talagrand's results on behavior of random processes.
The first one characterizes suprema of Gaussian and exponential processes in 
terms of the $\gamma_q$ functionals.

For a metric space $(E, \rho)$ and $q > 0$ we define the 
$\gam _q$ functional as 
\begin{displaymath}
\gamma_q (E,\rho) = \inf_{(A_s)_{s=0}^\infty} \sup_{x \in E}
\sum_{s=1}^\infty 2^{s/q}\ {\rm dist}(x,A_s),
\end{displaymath}
where the infimum is taken over all sequences $(A_s)_{s=0}^\infty$
of subsets of $E$, such that $|A_0| = 1$ and $|A_s| \le 2^{2^s}$
for $s \ge 1$.

The following theorem combines Theorems 2.1.1 and 5.2.7 in \cite{Tal}. 

\begin{theorem} 
\label{mm}
Let $T\subset \R^d$ and $\rho _q$ denote the $\ell _q$ metric. Then
$$
    \E \sup _{z\in T} \sum _{i=1}^d z_i g_i \approx \gam _2(T, \rho _2)
\quad \mbox{ and } \quad
 \E \sup _{z\in T} \sum _{i=1}^d z_i E_i \approx \gam _2(T, \rho _2)  +
  \gam _1(T, \rho _{\infty}).
$$
\end{theorem}

\medskip

We will also use Talagrand's result on the deviation of supremum of 
exponential processes from their averages. It follows by Talagrand's two 
level concentration for product exponential measure 
(\cite{Taltwo}).

\begin{theorem} \label{conexp}
Let $T$ be a compact subset of  $\R^d$. Then for any $t \ge 0$,  
\begin{displaymath}
    \P\left(\sup_{z\in T}\left|\sum_{i=1}^d z_i E_i\r| \ge \E 
    \sup_{z\in T}\bigg|\sum_{i=1}^d z_i E_i\bigg| + t\r) \le 
    \exp\left(-c\min\left\{\frac{t^2}{a^2},\frac{t}{b}\r\}\r),
\end{displaymath}
where $a = \sup_{z\in T} |z|$, $b = \sup_{z\in T}\|z\|_\infty$.
\end{theorem}

\section{Chevet type inequality}
\label{chevetineq}

\begin{theorem}\label{Chevet}
Let $\Gamma$ be an isotropic log-concave unconditional random 
$n\times N$ matrix. Let $K\subset \R^N$, $L\subset \R^n$ be 
origin symmetric convex bodies. Then 
\begin{align*}
 &\E \|\Gam : K \to L \|
\\
  &\leq
    C \left(\|\Id : K\to B_2^N \|  \cdot \E\left\|
   \sum_{i=1}^n E_i e_i \r\| _{L}  +  \|\Id : B_2^n \to L \|  \cdot \E\left\|
   \sum_{i=1}^N E_i e_i \r\| _{K^{\circ}} \r) .
\end{align*}
\end{theorem}


\medskip

\noindent
{\bf Example.}
One of the most important examples of matrices satisfying the hypothesis  
of Theorem~\ref{Chevet} are matrices whose rows (or columns) are independent 
isotropic log-concave unconditional random vectors. Indeed, it is easy to see 
that if $X$, $Y$ are independent isotropic log-concave  random vectors 
then so is $(X, Y)$. If $X$, $Y$ are in addition unconditional
then clearly $(X, Y)$ is unconditional. Therefore, if rows (or columns) of 
a matrix $\Gamma$ are independent isotropic log-concave random vectors then 
$\Gamma$ is isotropic log-concave. If rows (resp. columns) are in addition 
unconditional, then so is $\Gamma$. We will use it in Section~\ref{RIP}.

\bigskip

\noindent
{\bf Remarks. 1. }
In fact in the Gaussian case the equivalence holds in the 
Chevet inequality. However, in the log-concave case one cannot hope 
for the reverse inequality even in the case of exponential matrix and 
unconditional convex bodies $K$, $L$. Indeed, consider the matrix 
$\Gamma =(E_{ij})$ as an operator $\ell _1^N \to \ell _1^n$, 
i.e. $K=B_1^N$, $L=B_1^n$. By Lemma~\ref{lonenorm} 
$$
  \Ex\| \Gam \ :\ \ell _1^N \to \ell _1^n \| \approx n + \ln N . 
$$
On the other hand, the right hand side term in Theorem~\ref{Chevet} is 
$$
   C \left( \E \sum_{i=1}^n |E_i |  +  \sqrt{n}\  \E\max_{j\leq N}
   |E_i| \r) \approx n + \sqrt{n}\ \ln (2N). 
$$
Thus, if $N\geq e^n$ then the ratio between the right hand side and the 
left hand side is of the order $\sqrt{n}$. 
\\
{\bf 2.} The following weak form of a reverse inequality holds for the 
exponential matrix $\Gamma =(E_{ij})_{i\leq n,j\leq N}$:
\[
\E \|\Gam : K \to L \|
\geq \frac{1}{2}
\left(\max_{i\leq N}\|e_i\|_{K^{\circ}}  \cdot \E\left\|
   \sum_{i=1}^n E_i e_i \r\| _{L}  
+  
\max_{i\leq n}\|e_i\|_{L}  \cdot \E\left\|
   \sum_{i=1}^N E_i e_i \r\| _{K^{\circ}} \r) .
\]
Indeed,  fix $1\leq \ell\leq N$ and take $x\in K$ such that 
$\|e_\ell\|_{K^{\circ}}=|\langle e_\ell,x\rangle|=|x_\ell|$. 
Then
\begin{align*}
   \E \|\Gam : K \to L \|&\geq \Ex\|\Gam x\|_L=\Ex\Big\|\sum_{
    i\leq n, j\leq N} E_{ij}x_je_i\Big\|_L\geq \Ex\Big\|\sum_{
    i\leq n}E_{i\ell}x_\ell e_i\Big\|_L
\\
&=|x_\ell|\ \Ex\Big\|\sum_{i\leq n}E_ie_i\Big\|_L
=\|e_\ell\|_{K^{\circ}}\ \Ex\Big\|\sum_{i\leq n}E_ie_i\Big\|_L.
\end{align*}
This shows that 
$$
  \E \|\Gam : K \to L \|\geq \max_{i\leq N}\|e_i\|_{K^{\circ}}\ 
  \Ex\left\|\sum_{i\leq n}E_ie_i\r\|_L
$$
and by duality we have 
$$
  \E \|\Gam : K \to L \|= \E  \|\Gam^T : L^{\circ} \to K^{\circ} \| \geq 
  \max_{i\leq n} \|e_i\|_{L}\ \Ex\left\|\sum_{i\leq N}E_ie_i\r\|_{K^{\circ}} .
$$
{\bf 3. } As in Theorem~\ref{Latala}, the condition ``$\Gamma$ is unconditional" 
cannot be omitted  in Theorem~\ref{Chevet}. We show an example proving that 
in Section~\ref{example}. 

\medskip

\noindent
{\bf Proof of Theorem~\ref{Chevet}. }  
 First note that considering the matrix $\Gamma$ as a vector 
in $\R^{nN}$ and applying Theorem \ref{Latala}, we obtain that 
it is enough to prove Theorem~\ref{Chevet} for the case 
of the exponential matrix.

{}From now we assume that  $\Gam =(E_{ij})$. 
Denote as before
$
   T = K\otimes L^{\circ} = \{ x\otimes y \ \colon \ x\in K,\ y\in L^{\circ}  \}.
$
Then by Theorem~\ref{mm}
$$
  \E  \|\Gam : K \to L \| = \E  \sup _{x\in K} \sup _{y\in L^{\circ}}
  \sum _{i, j} E_{i j} x_j y_i = \E \sup _T \la \Gam, x\otimes 
  y \ra \approx \gam _2(T, \rho _2)  +   \gam _1(T, \rho _{\infty}) 
$$
and 
$$
   \E\left\| \sum_{i=1}^n E_i e_i \r\| _{L}  \approx
    \gam _2(L^{\circ}, \rho _2)  +   \gam _1(L^{\circ}, \rho _{\infty}) , 
$$
$$
   \E\left\| \sum_{i=1}^N E_i e_i \r\| _{K^{\circ}}
   \approx
   \gam _2(K, \rho _2)  +   \gam _1(K, \rho _{\infty}).
$$
Thus it is enough to show that
\begin{equation}\label{chinone}
    \gam _2(T, \rho _2)
   \leq  C\left( R(K) \gam _2(L^{\circ}, \rho _2)  +
    R(L^{\circ}) \gam _2(K, \rho _2)  \r)
\end{equation}
 and
\begin{equation}\label{chintwo}
   \gam _1(T, \rho _{\infty}) \leq  C\left(  R(K) \gam _1(L^{\circ}, \rho _{\infty})
     + R(L^{\circ})    \gam _1(K, \rho _{\infty}) \r).
\end{equation}

Inequality (\ref{chinone}) is the Chevet inequality for the Gaussian case.
Indeed by Theorem \ref{mm}
$$
    \gam _2(T, \rho _2)
     \approx \E \sup _{z\in T}
    \sum _{i, j} z_{ij} g_{ij} = \E \| (g_{ij})  \ : K \to L \| 
$$
and
$$
R(K)   \gam _2(L^{\circ}, \rho _2)  +
R(L^{\circ})  \gam _2(K, \rho _2)
\approx  R(K)  \E \sup _{z\in L^{\circ}} \sum _{i=1}^n z_i g_i
    + R(L^{\circ})     \E \sup _{z\in K} \sum _{i=1}^N z_i g_i .
$$
In fact we could prove (\ref{chinone}) without the use of the Chevet inequality, 
but by the chaining argument similar to the one used for the proof of 
(\ref{chintwo}) below (cf. also \cite{MT}).

It remains to prove inequality (\ref{chintwo}). 

Let $A_s \subset K$ and $B_s \subset L^{\circ}$, $s\geq0$,  be admissible 
sequences of sets (i.e., with $|A_0|=|B_0| = 1$, $|A_s|,|B_s| \le 2^{2^s}$ 
for $s \ge 1$). 
Define an admissible sequence $(C_s)_{s\geq 0}$ by $C_0 =\{0\}$ and 
$$
  C_s = A_{s-1}\otimes B_{s-1} \subset K\otimes L^{\circ}, \quad s\geq 1.
$$
Note that for all $x,\tilde{x} \in K$ and for all $y,\tilde{y}\in L^{\circ}$
one has 
\begin{align*}
   \|x\otimes y - \tilde x \otimes \tilde y \|_{\infty }
&\le
   \| x \|_{\infty } \cdot    \|y - \tilde y \|_{\infty }
     + \| \tilde y \|_{\infty }\cdot  \|x - \tilde x \|_{\infty }
\\
& \leq
  R(K) \|y - \tilde y \|_{\infty} +
  R(L^{\circ}) \|x - \tilde x \|_{\infty }.
\end{align*}
Therefore
\begin{align*}
\gamma_1(K\otimes L^{\circ}, \rho_{\infty}) 
\le&
  \sup _{x\otimes y \in K\otimes L^{\circ}}
  \sum_{s=0}^\infty 2^s\dist(x\otimes y,C_s)
\\  
\le&  R(K) \sup _{y \in L^{\circ}} \left( \|y\| _{\infty} + 
  \sum_{s=1}^\infty 2^s \dist(y, B_{s-1}) \r) 
\\
 &+R(L^{\circ}) \sup _{x \in K}  \left( \|x\| _{\infty} + 
  \sum_{s=1}^\infty 2^s \dist(x, A_{s-1})\r) . 
\end{align*} 
Taking the infimum over all admissible sequences $(A_s)$ and $(B_s)$ we get
\begin{align*}
\gamma_1&(K\otimes L^{\circ}, \rho_{\infty}) 
\\
&\le
  R(K) \left(\mbox{diam} _\infty L^{\circ} + 2 \gamma_1(L^{\circ}, \rho_\infty)\r)
  +  R(L^{\circ})\left(\mbox{diam} _\infty K + 2 \gamma_1(K, \rho_\infty)\r)
\\
 &\le 4 R(K) \gamma_1(L^{\circ}, \rho_\infty) + 4 R(L^{\circ}) \gamma_1(K, \rho_\infty) ,
\end{align*}
where in the last inequality we used the fact  that the diameter is clearly dominated 
by doubled $\gam _1$ functional. 
\qed

\begin{cor}\label{probest} Let $\Gamma$, $K$, $L$ be as in Theorem~\ref{Chevet}. 
Then for every $t>0$,  
$$
   \|\Gam : K \to L \|
  \leq  C \left( R(K)  \cdot \E\left\|
   \sum_{i=1}^n E_i  e_i\r\| _{L}  +  R(L^{\circ})  \cdot \E\left\|
   \sum_{i=1}^N E_i e_i\r\| _{K^{\circ}}   + t \r)
$$
with probability at least
$$
  1-\exp\left(-c\min\left\{\frac{t^2}{\sig^2},\frac{t}{\sig'}\r\}\r)\geq 
  1 - \exp\left(-c\min\left\{\frac{t^2}{\sig^2},\frac{t}{\sig}\r\}\r), 
$$
where $\sig = R(K)R(L^{\circ})$ and $\sig'=\sup_{x\in K}\|x\|_{\infty}\sup_{y\in 
L^{\circ}}\|y\|_{\infty}$. 
\end{cor}


\noindent
{\bf Proof.} 
As in the proof of Theorem \ref{Chevet} it is enough to consider the case  
$\Gam = (E _{ij})$. Moreover it suffices to show that 
\[
    \Pr(\|\Gam : K \to L \|\geq \Ex \|\Gam : K \to L \|+t)\leq 
    \exp\left(-c\min\left\{\frac{t^2}{\sig^2},\frac{t}{\sig'}\r\}\r). 
\]

To obtain the above estimate we use Theorem~\ref{conexp}. Recall that
$\|\Gam : K \to L\| =  \sup_T\langle \Gam,x\otimes y\rangle$, where 
$T=K\otimes L^{\circ}$. Thus we can easily compute parameters $a$ and 
$b$ in Theorem~\ref{conexp}: 
$$
  a = \sup _T  |x\otimes y | = \sup _{x\in K, \ y\in L^{\circ}}  |x|\cdot |y|
  = \sig 
$$
and
$$
 b = \sup _T \|x\otimes y \|_{\infty}=\sup _{x\in K, \ y\in L^{\circ}} 
 \|x\|_{\infty}\cdot \|y\|_{\infty} = \sig' .
$$ 
\qed

\section{Norms of submatrices and RIP}
\label{RIP}

Here we estimate the norms of submatrices of an isotropic unconditional 
log-concave random $n\times N$ matrix $\Gamma$.


 Recall that for subsets $ J\subset \{1,\ldots,n\}$ and $ I \subset
\{1,\ldots,N\}$,  $\Gamma(J, I)$  denotes the submatrix of $\Gamma$
consisting of the rows indexed by elements from $J$ and the columns
indexed by elements from $I$. Recall also that for $k\le n$ and $m\le N$,   
$\Gamma_{k,m}$ is defined by 
\begin{equation}
\label{akmcorr}
  \Gamma_{k,m} = 
  \sup  \|\Gamma (J, I)\ :\ \ell_2^m \to \ell_2^k\| ,
\end{equation}
where the supremum is taken over all subsets
$  J\subset \{1,\ldots,n\}$ and $ I \subset \{1,\ldots,N\}$  with
cardinalities  $ |J| = k, |I|= m$.
That is, $\Gamma_{k,m}$ is  the maximal operator norm 
of a submatrix of $\Gam$ with $k$ rows and $m$ columns.

We also denote the set of  $\ell$-sparse unit vectors on $\R^d$ 
by $U_{\ell}$ (or $U_{\ell}(d)$, when we want to emphasize the dimension 
of the underlying space) and its convex hull by $\tilde U_{\ell}$, i.e.
$$
  U_{\ell} = U_{\ell}(d) = \{x \in \R^{d} \colon |\supp x| \le \ell \, \, 
  \mbox{ and } \, \, |x|=1\},\quad
  \mbox{and}\quad  \tilde {U_{\ell}}=\conv(U_{\ell}).
$$
Thus 
$$
  \Gam_{k,m} = \left\| \Gam \ : \ \tilde U_m (N) \to (U_k (n))^{\circ} \r\| .
$$
Note that $(U_k (n))^{\circ} = (\tilde U_k (n))^{\circ}$. Below $U_{\ell}^{\circ}$ 
means $(U_{\ell})^{\circ}$.

\bigskip

\noindent {\bf Remark. }
For matrices with $N$ independent log-concave columns and $k=n$ the 
sharp estimates for $\Gamma_{n,m}$ were obtained in \cite{ALPT}.

\medskip

To treat the general case we will need the following simple lemma.

\begin{lemma}
\label{estUm}
For any $1\leq \ell\leq n$ we have
\[
 \E\left\|   \sum_{i=1}^n E_i e_i \r\| _{U_{\ell}^{\circ}} \approx 
 \sqrt{\ell}\ln\frac{3n}{\ell}.
\]
\end{lemma}


\noindent
{\bf Proof.} 
By Borell's lemma (\cite{Bo}) we have
\[
\left(\E\left\|   \sum_{i=1}^n E_i e_i \r\| _{U_{\ell}^{\circ}}\right)^2\approx
\E\left\|   \sum_{i=1}^n E_i e_i \r\|^2 _{U_{\ell}^{\circ}}
=\Ex\sup_{I\subset\{1.\ldots,n\}\atop |I|=\ell}\sum_{i\in I}E_i^2
=\sum_{i=1}^{\ell}\Ex |E_i^*|^2,
\]
where $E_1^*,\ldots,E_n^*$ denotes the nonincreasing rearrangement of 
$|E_1|,\ldots,|E_n|$. We conclude the proof by the standard well known 
estimate $\Ex |E_i^*|^2\approx (\ln (3n/i))^2$.
\qed

\medskip

Now observe that $\Gamma$ satisfies the hypothesis of Theorem~\ref{Chevet} 
and that $\tilde U_{\ell} \subset B_2^n$, so $R(\tilde U_{\ell})=1$. Thus  
Theorem~\ref{Chevet} implies 
$$
  \E   \Gamma_{k,m} \leq
    C \left(    \E\left\| \sum_{i=1}^N E_i e_i \r\| _{U_m^{\circ}}  +
   \E\left\|   \sum_{i=1}^n E_i e_i \r\| _{U_k^{\circ}} \r),
$$
which together with Lemma~\ref{estUm} and Corollary~\ref{probest} 
implies the following theorem. 

\begin{theorem}\label{subm} There are absolute positive constants 
$C$ and $c$ such that the following holds. 
Let $m\leq N$ and $k\leq n$. 
Let $\Gam$ be an isotropic unconditional 
log-concave random $n\times N$ matrix.  Then 
$$
   \E   \Gamma_{k,m} \leq
    C \left(    \sqrt{m}\ \ln\frac{3N}{m}  +
    \sqrt{k}\ \ln\frac{3n}{k} \r) .
$$
 Moreover, for every $t>0$,  
$$
   \Gamma_{k,m}  \leq  C \left(\sqrt{m}\ \ln\frac{3N}{m}  +
   \sqrt{k}\ \ln\frac{3n}{k} + t \r) 
$$
with probability at least
$$
   1 - \exp\left(-c\min\left\{ t, t^2 \r\}\r). 
$$
\end{theorem}

\medskip

\noindent
{\bf Remarks. 1.\ } In the case when $\Gam=(E_{ij})$ we have
\begin{align*}
  \E   \Gamma_{k,m} &\geq \max\Big\{\E\Big\|\sum_{i=1}^N E_i 
  e_i\Big\|_{U_m^{\circ}}, \E\Big\|\sum_{i=1}^n E_i e_i
  \Big\|_{U_k^{\circ}}\Big\}
\\
&\geq
\frac{1}{C} \left(    \sqrt{m}\ \ln\frac{3N}{m}  +
    \sqrt{k}\ \ln\frac{3n}{k} \r) .
\end{align*}

\noindent
{\bf 2. } Theorem~\ref{subm} can be proved directly 
(i.e. without Chevet inequality) using a chaining argument in 
the spirit of \cite{ALPT}. We provide the details in the last 
section. Similar estimates (with worse probability) were recently 
independently obtained in \cite{MP}.

\bigskip

We now estimate the restricted isometry constant (RIC) of a 
random matrix $\Gamma$ with independent unconditional isotropic
log-concave rows. As was mentioned in the example following 
Theorem~\ref{Chevet} such $\Gamma$ is unconditional isotropic
log-concave. Recall that the RIC of order $m$ is the smallest number 
$\delta = \delta _m(\Gam)$ such that
$$
   (1-\delta ) |x|^2 \leq |\Gam x|^2 \leq  (1+\delta ) |x|^2.
$$
for every $x\in U_m$.

%
%

The following theorem is an ``unconditional" counterpart of Theorem~6.4 
from \cite{ALLPT3} (see also Theorem~7 in \cite{ALLPT4}). Its proof repeats 
the lines of the corresponding proof in \cite{ALLPT3}. The result is sharp 
up to absolute constants.

\begin{theorem} \label{rip}  Let $0<\theta < 1$. 
Let $\Gam$ be an $n\times N$ random matrix, whose rows are independent 
unconditional isotropic log-concave vectors in $\R^N$. 
Then $\delta _m (\Gam/\sqrt n) \leq \theta$ with probability at least 
$$
   1-  \exp\left( - c\ \frac{ \theta ^2 n}{\ln ^2 n} \right) - 
     2 \exp{\left(- c\ \sqrt{m} \ln \frac{3 N }{m}\right)} ,
$$
provided that either 
\\
(i) \ $N\leq n$ and  
$$
   m\approx \min\left\{N, \ \frac{\theta ^2  n }{\ln ^3 (3/\theta) } \r\}
$$
or \\
(ii) \ $N\geq n$ and  
$$
  m \leq c \ \frac{\theta n}{\ln (3N/(\theta n))} \
    \min\left\{  \frac{1}{\ln (3N/(\theta n))},
    \frac{\theta }{\ln ^2 (3/\theta )}  \r\} ,  
$$
where $c>0$ is an absolute constant.  
\end{theorem}


\medskip

\noindent
{\bf Remarks. 1.} The condition on $m$ in (ii) can be written as follows 
$$
\mbox{if}\quad
\theta \geq \frac{\ln ^2 \ln (3N/ n)}{\ln (3N/ n)}
\quad \quad
\mbox{ then }
\quad \quad
m \leq c \ \frac{\theta n}{\ln ^2 (3N/(\theta n))} , 
\quad \quad\quad \quad\quad \quad \quad
$$
$$
\mbox{if}\quad  \theta \leq \frac{\ln ^2 \ln (3N/n)}{\ln (3N/n)}
\quad \quad
\mbox{ then }
\quad \quad
  m \leq c \ \frac{\theta ^2}{\ln^2 (3/\theta)} \
   \frac{n}{\ln (3N/ (\theta n))}.\quad \quad \quad \quad\, \,
$$
{\bf 2. } Precisely the proof of Theorem~6.4 in \cite{ALLPT3} 
(with estimates from our Theorem~\ref{subm}) gives that if 
$$
   b_m :=  m \left( \ln \frac{3 N }{m}\r)^2 \le c \theta n 
$$
and 
$$
  m \ln \frac{3 N}{m} \ln ^2 \frac{n}{b_m} \le c \theta ^2 n 
$$
then $ \delta _m (\Gam/\sqrt n) \leq \theta$ 
with probability at least 
$$
   1-  \exp\left( - c\ \frac{ \theta ^2 n}{\ln ^2 (n/b_m)} \right) - 
     2 \exp{\left(- c\ \sqrt{m} \ \ln \frac{3 N}{m}\right)} .
$$

\section{An example}
\label{example}

In this section we prove that the condition ``$X$ is unconditional" 
cannot be omitted in Theorems~\ref{Latala} and \ref{Chevet}. Namely, first 
we construct an example of isotropic log-concave non-unconditional 
$d$-dimensional random vector $X$ and a norm $\no$ on $\R^d$, which 
fails to satisfy the conclusion of Theorem~\ref{Latala}. Then we consider 
the matrix consisting of one column $X$ as an operator from $(\R, |\cdot|)$   
to $(\R^d, \no)$ and show that 
it does not satisfy the Chevet type inequality. 
The idea of the construction of $X$ is rather simple -- we start with 
a matrix with i.i.d. exponential entries and rotate its columns by a 
``random'' rotation. Considering the matrix as a vector with operator 
norm $\ell_1\to \ell _1$ we prove the result.

\begin{theorem} \label{latexam} 
Let $d\geq 1$ and $Y = (E_1, \ldots, E_d)$.  
There exists an isotropic log-concave random vector $X$ in $\R^d$ and a
norm $\|\cdot\|$ such that
\begin{equation}\label{lowbound}
 \E\|X\|  \ge c\ \sqrt{\ln d}\, \, \E\|Y\|, 
\end{equation}
where $c>0$ is an absolute constant. 
 Moreover, the $d\times 1$ matrix $B$, whose the only column is $X$, satisfies 
$$
 \E \| B : [-1, 1] \to L \|
  \geq
    c\ \sqrt{\ln d}\,  \left( \E\left\| \sum_{i=1}^d E_i e_i \r\| _{L} 
   + \|\Id : B_2^d \to L \| \r) , 
$$
where $L$ is the unit ball of $\no$. 
\end{theorem}

\medskip

 


\noindent
{\bf Proof. } Let $n, N$ be integers such that $d=nN$. 
Consider an $n\times N$ matrix $\Gam =(E_{ij})$. Denote its columns by 
$X_1, \ldots, X_N$, so  that $\Gam = [X_1, \ldots, X_N]$. As before, 
we consider $\Gam$ as a $d$-dimensional vector.  
Given $U\in O(n)$ rotate the columns of $\Gam$ by $U$:
$$
  A=A(U)=U \Gam = [UX_1, \ldots, UX_N]. 
$$ 
Then $A$ is a log-concave isotropic vector in $\R^d$. Below we show that 
if $N=\lfloor e^{cn} \rfloor$ for some absolute constant $c>0$ then there 
exists $U_0\in O(n)$ such that 
\begin{equation}\label{exalat}
  \E _{\Gam} \ \| A(U_0) \ :\ \ell _1^N \to \ell _1^n \|
  \geq c_1  \sqrt{\ln d}\  \
  \E _{\Gam} \ \| \Gam \ :\ \ell _1^N \to \ell _1^n \| .
\end{equation}
 This will prove the first part of the theorem, since it is clearly enough 
to consider only such $n, N, d$ by adjusting the constant in the main statement. 

To prove (\ref{exalat}) we estimate the average of $\|A(U)\|$ over $U\in O(n)$. 
For every $x$ in $\R^n$ we have
$$
  \P _{O(n)} \left(\left\{ \| U x \| _1 \geq c_2 \sqrt{n} \ \|x\|_2 \r\}\r)
 =\sigma_{n-1}(\{y\colon \|y\|_1\geq c_2\sqrt{n}\})
  \geq 1- \exp(- 2 c n),
$$
where $\sigma_{n-1}$ denotes the uniform distribution on $S^{n-1}$ and the last 
inequality follows by simple volumetric argument 
(or by concentration, see e.g. 2.3, 5.1 and 5.3 in \cite{MS}). 
Thus, if $N\leq e^{cn}$,
$$
  \P _{O(n)} \left(\left\{ \forall i\leq N \ : \| U X_i \| _1 \geq c_2  
\sqrt{n} \ \|X_i\|_2 \r\}\r)
  \geq 1- \exp(- c n) \geq \frac{1}{2}.
$$
Hence
$$
  \E _{O(n)} \ \max_{i\leq N} \| U X_i \| _1 \geq c_2 \sqrt{n} \
  \max_{i\leq N} \| X_i \|_2,
$$
which implies 
\begin{align*}
  \E _{\Gam} \ \E _{O(n)} \ \| A(U) \ :\ \ell _1^N \to \ell _1^n \| &\geq
  c_2 \sqrt{n} \ \ \E _{\Gam} \  \max_{i\leq N} \| X_i \|_2 \\
 &\geq
  c_2 \sqrt{n} \ \ \E _{\Gam} \  \max_{i\leq N} | E_{1, i} |
   \geq
  c_3 \sqrt{n} \ \ln N.
\end{align*} 
By Lemma~\ref{lonenorm} 
$$
  \E _{\Gam} \ \| \Gam \ :\ \ell _1^N \to \ell _1^n \| 
  \approx  n + \ln N  .
$$
Thus, taking $N=\lfloor e^{cn} \rfloor$,  
$$
  \frac{ \E _{O(n)} \ \E _{\Gam} \ \| A(U) \ :\ \ell _1^N \to \ell _1^n \|}{
   \E _{\Gam} \ \| \Gam \ :\ \ell _1^N \to \ell _1^n \|}
   \geq c_4 \frac{\sqrt{n} \ \ln N}{ n+\ln N}
 \geq c_5  \sqrt{\ln N}  \geq c_6  \sqrt{\ln d}. 
$$
Hence there exists $U_0\in O(n)$ satisfying (\ref{exalat}). 

Now we will prove the ``moreover" part of the theorem. Recall that $L$ is 
the unit ball of the norm $\|\cdot\|$ constructed above. The log-concave 
vector under consideration is $X=A(U_0)$ and the matrix which 
provides the counterexample to the Chevet type inequality is $B=[X]$. 
By the above calculations we have 
$$
  \E \| B : [-1, 1] \to L \| = \E \|X\| _L = 
  \E \| A(U_0) \ :\ \ell _1^N \to \ell _1^n \| 
   \geq  c\ (\ln d)^{3/2}  
$$
 and 
$$
  \E\left\| \sum_{i=1}^d E_i e_i \r\| _{L} = \E \| \Gam \ :\ \ell _1^N 
  \to \ell _1^n \| \approx  n + \ln N  \approx \ln d. 
$$
It is easy to check that for every $n\times N$ matrix $T=(t_{ij})$ 
one has 
$$
  \| T \ :\ \ell _1^N \to \ell _1^n \| = \max _{j\leq N} \sum_{i=1}^{n} 
  |t_{ij}|  \leq \sqrt{n} \left(\sum_{i=1}^{n}\sum_{i=1}^{n} |t_{ij}|^2
   \r)^{1/2} = \sqrt{n}  \ |T|, 
$$
where $\sqrt{n}$ is the best possible constant in the inequality. 
This shows that 
$$
  \|\Id : B_2^d \to L \| =\sqrt{n} \approx \sqrt{\ln d}.
$$
Thus 
$$
 \E\left\| \sum_{i=1}^d E_i e_i \r\| _{L} + \|\Id : B_2^d \to L \| 
 \approx  \ln d, 
$$
which completes the proof. 
\qed

\bigskip 


\noindent
{\bf Concluding remarks. 1.} 
The above example is  optimal in the sense that one can't expect 
better than $\sqrt{\ln d}$ dependence on dimension in (\ref{lowbound}).  
Indeed, let $Y=(E_1,\ldots,E_d)$. We show that for any $U\in O(d)$ and 
any norm  $\|\cdot \|$ on $\er^d$ one has  
\begin{equation}\label{contrex}
  \Ex\|UY\|\leq C\sqrt{\log(ed)}\ \Ex\|Y\|.
\end{equation} 
First it is known that $\Ex\|Y\|\leq C\sqrt{\log(ed)}\ \Ex\|G\|$, where 
$G=(g_1, \ldots, g_d)$. Now note that if $K$ is a unit ball of $\|\cdot\|_K$ 
then for every  $U\in O(d)$ one has $\|U x\|_K = \|x\| _{U^{-1}K}$ for 
every $x\in \R^d$. Therefore, for any $U\in O(d)$ we have 
$$
  \Ex\|U Y\|\leq C\sqrt{\log(ed)}\ \Ex\|U G\|= C \sqrt{\log(ed)} \ \Ex\|G\|
$$
(in the last equality we used that the distribution of $G$ is invariant under 
rotations). 
Finally note that by either Theorem~\ref{mm} or Theorem~\ref{Latala} the norm 
of an exponential random vector dominates the norm of the Gaussian one, i.e.
$\Ex\|G\|\leq C_1\,  \Ex\|Y\|$, which implies (\ref{contrex}).
\\
{\bf 2.} For any isotropic vector $X$ in $\R^d$ (not necessarily log-concave) 
and any origin symmetric convex body $K\subset \R^d$ we show that 
\begin{equation} \label{cdcon}
  \Ex\|X\|_K\leq Cd(K,B_2^d)\ \Ex\|Y\|_K, 
\end{equation}
where $Y=(E_1,\ldots,E_d)$ and $d(K,B_2^d)$ denotes the Banach-Mazur distance 
between $K$ and $B_2^d$. Since for every origin symmetric $K$ one has 
$d(K,B_2^d)\leq \sqrt{d}$ (see e.g. \cite{Tom}), the inequality (\ref{cdcon}) 
implies that for any norm $\|\cdot \|$ on $\R^d$ 
$$
  \Ex\|X\|\leq C\sqrt{d}\ \Ex\|Y\| .
$$
Now we prove (\ref{cdcon}). First, as in Remark~1, note that the norm 
of an exponential random vector dominates the norm of the Gaussian one. Thus 
it is enough to show that $\Ex\|X\|_K\leq Cd(K,B_2^d)\ \Ex\|G\|_K$, 
where $G$ is as in Remark~1. 
Let $\alpha=d(K,B_2^d)$ and ${\cal E}$ be an ellipsoid such that 
${\cal E}\subset K\subset \alpha {\cal E}$. Since this is only a matter of 
rotation of a coordinate system we may assume that ${\cal E}=\{x\in \R^d\colon 
\sum_{i=1}^d a_i^2x_i^2\leq 1\}$. Then by the isotropicity of $X$,
$$
  \Ex\|X\|_{K}\leq \Ex\|X\|_{{\cal E}}=\Ex\left(\sum_{i=1}^d a_i^2X_i^2\r)^{1/2}
  \leq \left(\sum_{i=1}^d a_i^2\r)^{1/2}\leq C\Ex\|G\|_{{\cal E}}\leq C\alpha
  \Ex\|G\|_K, 
$$ 
where we used comparison of the  first and second moments of the norm 
$\|G\|_{{\cal E}}$ of the Gaussian vector.

\section{A direct proof of Theorem~\ref{subm}}
\label{dirproof}

We present here a proof of Theorem~\ref{subm} not involving the 
Chevet type inequality and not relying on Theorem~\ref{mm}, 
but only on tail estimates for suprema of linear combinations of 
independent exponential variables given in Theorem~\ref{conexp}.

We need the following lemma, which is an immediate consequence of 
Theorem~\ref{conexp} (recall here that for a matrix $A=(a_{ij})$, 
$\|A\|_{\infty}$ denotes $\max _{i,j} |a_{ij}|$).

\begin{lemma}\label{individual_bound:lemma} 
For every $n\times N$ matrix $A=(a_{ij})$ and every $t\geq 0$ we have
\begin{displaymath}
 \PP\left(\left|\sum_{ij} E_{ij} a_{ij}\r| \ge t\r) \le 
 2\exp\left(-c\min\left(\frac{t^2}{|A|^2},\frac{t}{\|A\|_\infty}\r)\r), 
\end{displaymath}
where $c>0$ is an absolute constant. 
\end{lemma}

Indeed, since $\E |\sum_{ij} E_{ij} a_{ij}| \leq (\E |\sum_{ij} E_{ij} 
x_{ij}|^2)^{1/2} =|A|$, the above Lemma follows from Theorem~\ref{conexp} 
for $t\geq 2|A|$. For $t\leq 2|A|$ we can make the 
right hand side larger than 1 by the choice of $c$. 

\medskip

\noindent
{\bf Direct proof of Theorem~\ref{subm}. } 
As in the proof of Theorem~\ref{Chevet}, using Theorem~\ref{Latala}, we may 
assume that $\Gamma$ is the exponential matrix, i.e. $\Gamma = (E_{ij})$. 
Without loss of generality we assume that $k\geq m$ and that $k = 2^{r}-1$, 
$m = 2^{s}-1$ for some positive integers $r\geq s$. It is known (and easy to 
see by volumetric argument) that for any origin symmetric convex body 
$V\subset \R^d$ and any $\eps\leq 1$ there exist an $\eps$-net (with 
respect to the metric defined by  $V$) in $V$ of cardinality at most 
$(3/\eps)^d$. For $i = 0,1,\ldots,r-1$ let $\mathcal{M}_i$ be a 
$(2^{i}/(4k))$-net (with respect to the metric defined by $B_2^n \cap 
(2^{-i/2}B_\infty^n$)) in the set
\begin{displaymath}
\bigcup_{{I\subseteq \{1,\ldots,n\}}\atop{|I| \le 2^i}}
\R^I \cap B_2^n\cap (2^{-i/2}B_\infty^n)
\end{displaymath}
of cardinality not greater than
\begin{displaymath}
\binom{n}{2^i}\Big(\frac{12k}{2^i}\Big)^{2^i} \le
\exp\Big(C2^i\log\Big(\frac{2n}{2^i}\Big)\Big), 
\end{displaymath}
where $\R^I$ denotes the span of $\{e_i\}_{i\in I}$.
Similarly for $i = 0,1,\ldots,s-1$ let $\mathcal{N}_i$ be a
$(2^i/(4m))$-net in the set
\begin{displaymath}
\bigcup_{{I\subseteq \{1,\ldots,N\}}\atop{| I| \le 2^i}}
\R^I \cap B_2^N\cap (2^{-i/2}B_\infty^N)
\end{displaymath}
of cardinality at most
\begin{displaymath}
\binom{N}{2^i}\Big(\frac{12m}{2^i}\Big)^{2^i} \le
\exp\Big(C2^i\log\Big(\frac{2N}{2^i}\Big)\Big).
\end{displaymath}

Let now $\mathcal{M}$ be the set of vectors in $2B_2^n$ that can
be represented in the form $x = \sum_{i=0}^{r-1} x_i$, where $x_i
\in \mathcal{M}_i$ and have pairwise disjoint supports.
Analogously define $\mathcal{N}$ as the set of vectors $y =
\sum_{i=0}^{s-1}y_i\in 2B_2^N$, with $y_i \in \mathcal{N}_i$ and
pairwise disjoint supports. For $x \in \mathcal{M}$ and $i =
0,1,\ldots,r-1$ let $S_i x = x_0+\ldots+x_i$, where $x_i$ is the
appropriate vector from the above representation (this
representation needs not be unique, so for each vector $x$ we
choose one of them). Similarly, for $i  = 0,1,\ldots,s-1$ and $y
\in \mathcal{N}$ let $T_iy = y_0+\ldots+y_i$. For $i =
s,\ldots,r-1$ let $T_iy = y$. Additionally set $S_{-1} x = 0$,
$T_{-1}y = 0$. We thus have
\begin{displaymath}
 y \otimes x=  \sum_{i=0}^{r-1} (T_iy\otimes S_ix - T_{i-1}y\otimes S_{i-1}x)
\end{displaymath}
for $x \in \mathcal{M}, y \in \mathcal{N}$.

Since $x_i$'s and $y_i$'s have pairwise disjoint supports, viewing 
$(T_j y\otimes S_j x)$'s as sub-matrices of $y\otimes x$, it is easy to 
check that for every $j\geq i$  
\begin{equation}\label{twonorm}
 |T_j y\otimes S_j x  - T_{i-1}y\otimes S_{i-1}x| \le 4
\end{equation}
and 
\begin{equation}\label{infnorm}
 \|T_j y\otimes S_j x  - T_{i-1}y\otimes S_{i-1}x\|_{\infty} \le 2^{-i/2}. 
\end{equation} 
Thus, by Lemma \ref{individual_bound:lemma}, for any $x \in
\mathcal{M}$, $y \in \mathcal{N}$ and $t\ge 1$,
\begin{align}\label{individual_bound:eq} 
 \PP(|\langle \Gamma T_i y, S_i x\rangle - \langle \Gamma T_{i-1} y, 
 S_{i-1} x\rangle| \ge t) \le 2\exp(-c\min(t^2,2^{i/2} t)).
\end{align}

Moreover, for any $i\le s-1$, the cardinality of the set of
vectors of the form $T_i y\otimes S_i x  - T_{i-1}y\otimes
S_{i-1}x$, $x\in \mathcal{M}, y \in \mathcal{N}$ is at most 
$$
  \exp\Big(\sum_{j=0}^i\Big(C2^j\log\Big(\frac{2n}{2^j}\Big)+ 
  C 2^j\log\Big(\frac{2N}{2^j}\Big)\Big)\Big) 
  \le \exp\Big(\tilde{C}2^i\log\Big(\frac{2n}{2^i}\Big)+
  \tilde{C}2^i\log\Big(\frac{2N}{2^i}\Big)\Big).
$$

By (\ref{individual_bound:eq}) and the union bound we get that for
$i\le s-1$ and any $t \ge 1$, with probability at least
\begin{align*}
&1 - 2\exp\Big(-ct\Big(2^{i}\log(2n/2^{i}) +
2^{i}\log(2N/2^i)\Big)\Big),
\end{align*}
one has 
$$
  \max_{x\in\mathcal{M},y\in\mathcal{N}} |\langle \Gamma T_i y, S_i 
  x\rangle - \langle \Gamma T_{i-1} y, S_{i-1} x\rangle| 
  \le Ct\Big(2^{i/2}\log(2n/2^{i}) + 2^{i/2}\log(2N/2^i)\Big). 
$$
By integration this yields 
$$
  \E \max_{x\in\mathcal{M},y\in\mathcal{N}} |\langle \Gamma T_i y, 
  S_i x\rangle - \langle \Gamma T_{i-1} y, S_{i-1} x\rangle| 
  \le C\Big(2^{i/2}\log(2n/2^{i}) + 2^{i/2}\log(2N/2^i)\Big).
$$
Therefore 
\begin{align}\label{first_part}
\E\sup_{x\in\mathcal{M},y\in\mathcal{N}}|\langle \Gamma T_{s-1} y, 
 S_{s-1}x\rangle| &\le\sum_{i=0}^{s-1} \E \sup_{x\in\mathcal{M}, y\in
 \mathcal{N}}|\langle \Gamma T_i y, S_i x\rangle-\langle \Gamma T_{i-1} y, 
 S_{i-1} x\rangle| \nonumber\\
 &\le \sum_{i=0}^{s-1} C\Big(2^{i/2}\log(2n/2^{i}) +
 2^{i/2}\log(2N/2^i)\Big) \nonumber\\
&\le C_1 \Big(\sqrt{k}\log(2n/k) +
\sqrt{m}\log(2N/m)\Big).
\end{align}

On the other hand, for any $y \in \mathcal{N}$ and $i \ge s$, we
have by $T_{i-1}y = T_i y =y$. Thus by (\ref{individual_bound:eq})
and the fact that there are at most $\exp(C 2^i\log(2n/2^i))$
vectors of the form $S_i x - S_{i-1} x$ with $x \in \mathcal{M}$,
we get for $t\ge 1$,
\begin{displaymath}
\sup_{x\in \mathcal{M}} |\langle \Gamma T_i y, (S_i x -
S_{i-1}x)\rangle|\le C t 2^{i/2}\log(2n/2^i),
\end{displaymath}
with probability at least $1 - \exp(-ct2^i\log(2n/2^i))$.

This implies that for $s\le i \le r-1$,
\begin{displaymath}
\E \max_{x\in \mathcal{M}} |\langle \Gamma T_i y, S_i x\rangle - \langle
\Gamma  T_{i-1} y, S_{i-1} x\rangle|\le C 2^{i/2}\log(2n/2^i)
\end{displaymath}
and thus 
 \begin{align*}
   \E&\max_{x \in \mathcal{M}} |\langle \Gamma T_{r-1} y, S_{r-1}x\rangle - 
   \langle \Gamma T_{s-1} y, S_{s-1} x\rangle| \\
   &\le \sum_{i=s}^{r-1}\E\max_{x \in \mathcal{M}} |\langle \Gamma T_i y, 
   S_i x\rangle - \langle T_{i-1} \Gamma y, S_{i-1} x\rangle|\\
   &\le C \sum_{i=s}^{r-1} 2^{i/2}\log(2n/2^i) \le
   \tilde{C}\sqrt{k}\log(2n/k).
\end{align*}

Applying Theorem~\ref{conexp} together with (\ref{twonorm}) and 
(\ref{infnorm}) (with $j=r-1$ and $i=s$) 
we obtain that for any $y \in \mathcal{N}$ and $t \ge 1$,
\begin{displaymath}
 \max_{x \in \mathcal{M}} |\langle \Gamma T_{r-1}y, S_{r-1}x\rangle -
 \langle \Gamma  T_{s-1} y, S_{s-1} x\rangle| \le C\sqrt{k}\log(2n/k) +
 Ct2^{s/2}\log(2N/2^s),
\end{displaymath}
with probability at least
\begin{displaymath}
  1 - 2\exp(-\tilde{C}t 2^s\log(2N/2^s)),
\end{displaymath}
which by the union bound and integration by parts gives
\begin{align*}
\E&\max_{x\in \mathcal{M},y\in\mathcal{N}} |\langle T_{r-1}y,A S_{r-1}x\rangle 
- \langle T_{s-1} y,AS_{s-1} x\rangle| \\
&\le C\sqrt{k}\log(2n/k) + C2^{s/2}\log(2N/2^s) \le
\tilde{C}\Big(\sqrt{k}\log(2n/k)+\sqrt{m}\log(2N/m)\Big).
\end{align*}

Combining this inequality with (\ref{first_part}) we get
\begin{displaymath}
 \E\max_{x\in\mathcal{M},y\in\mathcal{N}} |\langle y, A x\rangle|
 \le C\Big(\sqrt{k}\log(2n/k)+\sqrt{m}\log(2N/m)\Big).
\end{displaymath}

Let us now notice that for arbitrary $x \in S^{n-1}$, $y \in
S^{n-1}$, with $|\supp x| \le k, |\supp y| \le m$, there exist
$\tilde{x}\in \mathcal{M}, \tilde{y} \in \mathcal{N}$, such that
$\supp \tilde{x} \subset \supp x$, $\supp \tilde{y} \subset \supp
y$ and
\begin{align*}
  |x - \tilde{x}|^2 \le \sum_{i=0}^{r-1} 2^{2i}/(16k^2) \le 1/8,
  \quad|y - \tilde{y}|^2 \le \sum_{i=0}^{s-1} 2^{2i}/(16m^2) \le 1/8.
\end{align*}

We have
\begin{displaymath}
 \langle \Gamma  y, x\rangle = \langle \Gamma  \tilde{y}, \tilde{x}\rangle +
 \langle \Gamma (y - \tilde{y}), x\rangle + \langle \Gamma \tilde{y}, x-\tilde{x} 
 \rangle.
\end{displaymath}
Taking into account that $\tilde{y} \in 2B_{2}^N$ and passing to suprema, we get
\begin{displaymath}
 \Gamma _{k,m} \le \max_{\tilde{x}\in\mathcal{M},\tilde{y}\in\mathcal{N}}
 \langle \Gamma  \tilde{y}, \tilde{x}\rangle + 3 \Gamma _{k,m} /8 
\end{displaymath}
and thus
\begin{displaymath}
 \E \Gamma _{k,m} \le 2\E
 \max_{\tilde{x}\in\mathcal{M},\tilde{y}\in\mathcal{N}} \langle
 \Gamma \tilde{y}, \tilde{x}\rangle\le C(\sqrt{k}\log(2n/k) + \sqrt{m}\log(2N/m)),
\end{displaymath}
which completes the proof of the first part of Theorem~\ref{subm}. 
The proof of the ``moreover" part is obtained using Theorem~\ref{conexp} in the same 
way as it was used  to obtain Corollary~\ref{probest} from Theorem~\ref{Chevet}.  
\qed

\medskip

\noindent
{\bf Remark.} We would like to notice that by adjusting the chaining 
argument presented above one can eliminate the use of the full strength 
of Theorem~\ref{conexp} and obtain a proof relying only on tail inequalities 
for linear combinations of independent exponential random variables 
(which follow from classical Bernstein inequalities). The modification 
involves splitting the proof into two cases depending on the comparison 
between $m\log(2N/m)$ and $k\log(2n/k)$.

\vspace{1cm}

\address

\end{document}